\documentstyle{amsppt}    

\loadbold
\loadmsam
\loadmsbm

\def\N{{\Bbb N}}

\def\R{{\Bbb R}}
\def\C{{\Bbb C}}

\def\P{{\Bbb P}}

\topmatter

\title
 Differential properties of  matrix orthogonal  polynomials
\endtitle
\leftheadtext{ Differential properties of  matrix orthogonal  polynomials}
\rightheadtext{M.J. Cantero, L. Moral, L. Vel\'azquez}

\author  
M.J. Cantero, L. Moral, L. Vel\'azquez
\endauthor

\affil
Departamento de Matem\'atica Aplicada,
Universidad de Zaragoza,
Spain
\endaffil

\date  
 December 2001
\enddate

\thanks
This work  was supported by Direcci\'on General de
Ense\~nanza Superior (DGES) of Spain under grant PB 98-1615. 
 \endthanks

\endtopmatter

\noindent {\bf Abstract}

  In this paper a general theory of semi-classical matrix orthogonal polynomials is developed.
We define the semi-classical linear functionals by means of a distributional equation $D(u A) = u
B,$ where $A$ and $B$ are matrix polynomials. Several characterizations for these semi-classical
functionals are given in terms of  the corresponding (left) matrix orthogonal polynomials
sequence. They involve  a   quasi-orthogonality property for their derivatives, 
a structure relation and
a second order differo-differential equation. Finally we illustrate the preceding results with 
  some non-trivial examples. 
 
\medskip

{\noindent \it Keywords and phrases:}    Orthogonal matrix polynomials, Semi-classical
functionals, Differential equation, Structure relation. 

{\noindent\it (2000) AMS Mathematics Subject Classification  }: 42C05 
\vskip 0.5cm
 
{\noindent \it Suggest running head}: Differential properties of matrix O P.
\vskip 0.25cm
{\noindent \it Corresponding author}: lmoral\@posta.unizar.es
 
$\;\;\;\;\;\;\;\;\;\;\;\;\;\;\;\;\;\;\;\;\;\;${\it Tel}: +34 976 76 11 41
 
$\;\;\;\;\;\;\;\;\;\;\;\;\;\;\;\;\;\;\;\;\;\;${\it Fax}: +34 976 76 11 25

\vfill\eject
\null\vskip 0.5cm

{\bf \S $\;\;\;$ 1 -  Introduction. }

\medskip 

  The study of   semi-classical orthogonal polynomials  in the scalar case (i.e., the study of
orthogonal polynomial whose associated functional satisfies a distributional equation $D(u\Phi) =
u\Psi$, where $\Phi,$ $\Psi$ are polynomials  with deg $\Psi \geq 1$)  was started by Shohat
([15]) in order to generalize the    properties of classical orthogonal polynomials.

  Among others, in ([2, 9]),  an approach to such polynomials taking into account the 
quasi-orthogonality of the derivatives of the polynomials sequence  is given. We also mention  the results
of Maroni ([11, 12])  where   an algebraic theory  of semi-classical orthogonal
polynomials  is presented.  The purpose  of
this theory is
 the characterization of semi-classical orthogonal polynomials by means of the quasi-orthogonality of
their derivatives, the structure relation and the differo-differential equation of  second
order.  

  In the last years, the study of matrix orthogonal polynomials attracted a
great  interest of the researchers, (see [4, 6, 10]).  As for their differential properties it is
known  the result of Dur\'an ([5]) who characterizes those matrix orthogonal polynomials  (O P)
satisfying a symmetric second order differential equation with polynomial coefficients. He proves
that they are diagonal (up  to a factor) with  classical  scalar O P in the diagonal. In
spite of the variety of applications of matrix polynomials ([4, 5, 6]), there are not too many
known families of semi-classical matrix O P  out of the diagonal case.

 One way to study many families of matrix orthogonal polynomials 
is to extend the analysis started by Dur\'an of differential properties of matrix orthogonal
polynomials. A natural way to do this is to generalize the theory of semi-classical  scalar O P to
the matricial case. In order to do that, we start with a Pearson type equation for   matrix
quasi-definite functionals.  We   obtain  their characterization in terms of a (in general,
non-symmetric) differo-differential equation for the related matrix orthogonal polynomials. In
the way of the proof,  we find another equivalences that generalize the scalar case. 

 This paper has been organized as  follows.
In   Section 3   we introduce  and study the concept of quasi-orthogonality   for
  matrix   polynomials.

 The mean  result of  Section 4 is the caracterization of semi-classical functionals in terms of
the quasi-orthogonality of the corresponding matrix   polynomials and their  derivatives.

 In Section 5 we  prove  for matrix O P with the above  quasi-orthogonality pro\-per\-ties, the
structure relation  and the   differo-differential equation, showing that they are equivalent to a
 Pearson-type equation. 

  Finally, in 
Section 6, we illustrate the preceding with some     examples  and we find  a way   
  to construct non-diagonalizable   semi-classical matrix functionals.

 %\vfill\eject
\vskip 0.5cm
{\bf \S$\;\;\;$ 2 - Basic tools. }

\medskip

We shall denote by $\P^{(m)}\;\,$ the $\C^{(m,m)}$-left-module
$$
\P^{(m)}=\left\{ \sum_{k=0}^n \alpha_k x^k\big | \alpha_k\in \C^{(m,m)}; n\in\N \right\}
$$
and by means of $\P^{{(m)}'}\;$   the  $\C^{(m,m)}$-right-mo\-du\-le   
Hom$\left(\P^{(m)},
\C^{(m,m)}\right).\,\,$  
  
\vskip 0.25cm

{\bf 2.1 Definition.} {\sl

(i)  The duality bracket is defined by  
 $$
\eqalign{ \bigl< \cdot\;,  \;\cdot \bigr> :   \P^{(m)} &\times \P^{{(m)}'} \to \C^{(m,m)} \cr
                                             \bigl (P&, u \bigr) \to \bigl< P, u \bigr> :=   
u\bigl(P\bigr) \cr}
$$

(ii) For $k\in \N$ and   $u\in \P^{{(m)}'}$ the linear functional $u x^k I\in  \P^{{(m)}'}$ is
defined by   
$$\bigl < P, u x^kI\bigr > :=    \bigl < x^k P, u  \bigr >,
$$
where $I$ denotes the $m\times m$ identity matrix.

  A linear extension gives
the right-product $\;u\cdot A \in \P^{{(m)}'}\;$ for $\;u\in \P^{{(m)}'},\;$  $\;A\in
\P^{(m)},\;$  with $\;\;\displaystyle A(x)=\sum_{k=0}^n\alpha_k x^k,\;$
in the following way:
$$
\bigl < P, u A \bigr > =  \sum_{k=0}^n \bigl < x^k P, u \bigr > \alpha_k.
$$

(iii)  The inner product  
$$\eqalign{ \bigl< \cdot\;, \;\cdot \bigr> :   \P^{(m)}&\times \P^{(m)} \to \C^{(m,m)} \cr
                                            \bigl (P&, Q \bigr) \to \bigl< P, Q \bigr>_u :=
\bigl<P, u Q^*\bigr> \cr}
$$
is defined for every  $u\in \P^{{(m)}'},$  where $Q^*$   denotes the trasposed conjugated of
$Q$.}

\vskip 0.25cm

{\bf  Remark.} 

  The duality bracket verifies the following linear properties: 

$$
\bigl< \alpha_1 P_1 + \alpha_2 P_2,\;  \beta_1 u_1 + \beta_2 u_2 \bigr>\; = 
$$
$$
=\;  \alpha_1 \bigl< P_1, u_1 \bigr> \beta_1 
+ \alpha_1 \bigl< P_1, u_2 \bigr> \beta_2  
+ \alpha_2 \bigl< P_2, u_1 \bigr> \beta_1  
+ \alpha_2 \bigl< P_2, u_2 \bigr> \beta_2,  
$$ 
for  all
 $\alpha_1, \alpha_2,\beta_1,\beta_2\in
\C^{(m,m)},\;$  $P_1, P_2 \in \P^{(m)}$   and $u_1, u_2\in \P^{{(m)}'},$  while the inner product is
sesquilinear:

(i) $\bigl< \alpha_1 P_1 + \alpha_2 P_2, \beta_1 Q_1 + \beta_2 Q_2\bigr >_u =  \alpha_1 \bigl< P_1, Q_1 \bigr>_u \beta_1^* + 
   + \alpha_1 \bigl< P_1, Q_2 \bigr>_u \beta_2^* +$

\centerline{ $ + \alpha_2 \bigl< P_2, Q_1 \bigr>_u \beta_1^* + \alpha_2 \bigl< P_2, Q_2 \bigr>_u \beta_2^*, $}

{\noindent for}  all
 $\alpha_1, \alpha_2,\beta_1,\beta_2\in
\C^{(m,m)},\;$  $P_1, P_2, Q_1, Q_2 \in \P^{(m)};$

(ii) $ \bigl< P, Q \bigr >_u^* =  \bigl< Q, P \bigr >_u$ 

{\noindent for}  all $ P, Q  \in \P^{(m)}.$

\vskip 0.25cm

{\bf 2.2 Definition.} {\sl We denote by   $ C_k :=\bigl< x^kI, u \bigr>\;\; $   the $k$-th 
moment with respect to  $u\in \P^{{(m)}'}$. Given  $u\in \P^{{(m)}'}$ with moments
$\left(C_k\right)_{k\in \N}, $ we 
say that $u$  is  quasi-definite 
  (non-singular)    if
det$\bigl[\bigl(C_{k+j}\bigr)_{k,j=0}^n\bigr]\not=0,\;$ $\forall n\geq 0,$  were 
$\left(C_{k+j}\right)_{k,j=0}^n$ is the Hankel-block matrix
$$  
\pmatrix 
     C_0 & C_1  &\cdots &C_n   \\ 
     C_1 &C_2&\cdots &C_{n+1}  \\ 
     \cdots&\cdots&\cdots&\cdots  \\
     C_n&C_{n+1}&\cdots&C_{2n} 
\endpmatrix
$$}
 
\vskip 0.25cm

{\bf Remark.} Given the sequence $\;\;\bigl(C_k\bigr)_{k=0}^\infty \subset
\C^{(m,m)}\;\;$, there exists a unique $u\in \P^{{(m)}'}$ such that $\bigl< x^kI, u \bigr > = C_k$.
This  establishes  an isomorphism between $\P^{{(m)}'}$ and the formal series with
coefficients in    
$\;\C^{(m,m)},\;\;\;$ $\;\;\displaystyle \sum_{k=0}^\infty C_k x^k$.

\vskip 0.25cm

{\bf 2.3 Definition.} {\sl  Let $u\in \P^{{(m)}'}.$  We say that $u$ is hermitian if $C_k^* =
C_k$,  for all $k\geq 0$.}

{\bf 2.4 Theorem.} {\sl Let  $u\in \P^{{(m)}'}$ be quasi-definite. Then, there exists a unique
 (up to left non-singular matrix factors)
  sequence of left orthogonal polynomials $\bigl(P_n\bigr)_{n\geq 0}$ with respect to $u,$ that is

\item{(1)} $ P_n\in\P^{(m)},$ dg$P_n = n.$

\item{(2)} The leading coefficient of $P_n$ is non-singular.

\item{(3)} $\bigl<P_n, P_m\bigr>_u = K_n\delta_{nm}$, where $K_n$ is non-singular.

 This sequence verifies a recurrence relation
$$
\eqalign{&x P_n(x) = \alpha_n   P_{n+1}(x) + \beta_n P_n(x) + \gamma_n P_{n-1}(x), \qquad n\geq 0,
\cr
         &P_{-1}(x) = \theta,  \cr} \eqno (1.1)
$$
%$$
%\eqalign{&P_{n+1}(x) = (\alpha_n x + \beta_n) P_n(x) + \gamma_n P_{n-1}(x), \qquad n\geq 0, \cr
%         &P_{-1}(x) = \theta,  \cr} \eqno (1.1)
%$$

 {\noindent  with} $\alpha_n$, $\beta_n$, $\gamma_n\in \C^{(m,m)},$    $\alpha_n$, $\gamma_n$
  non-singular and we denote by $\theta,$ the zero matrix.}

{ \bf{\sl  Proof.} } See [4, 14]. $\diamond\diamond$

\vskip 0.25cm

{\bf Remarks.}

\item{(i)} If we consider   the structure of right-modulus  on    $
\P^{(m)}$ (so,  left-modulus on   $\P^{{(m)}'}$), we define in  a similar way  a 
  sequence  of right orthogonal polynomials.

\item{(ii)} The inner product $\bigl<P_n, x^n I\bigr>_u$ is non-singular  for every quasi-definite
$u \in
\P^{{(m)}'}.$

\vskip 0.75cm

{\bf \S $\;\;\;$ 3 -  Quasi-orthogonality. }

\vskip 0.25cm

  Like in the scalar case, the semi-classical character will be closely related to the idea of
quasi-orthogonality. But the scalar quasi-orthogonality, when ge\-ne\-ra\-li\-zed to the matricial
case, splits in two different concepts. Both of them will play an important role in the
characterization of semi-classical matrix functionals.

\vskip 0.25cm

{\bf 3.1 Definition.} {\sl Let $\bigl(Q_n\bigr)_{n\geq 0}$ be a sequence of matrix  polynomials  
such that the leading coefficient of $Q_n$ is  non-singular 
and dg $Q_n=n$.

{\bf 1.-} We will say  that $\bigl(Q_n\bigr)_{n\geq 0}$ is a sequence of  left quasi-orthogonal 
matrix polynomials of order
$r$ with respect to $v \in \P^{{(m)}'}\setminus\{0\}$ if 
\item{(a)} $\bigl<x^k Q_n, v\bigr> = \theta, $\qquad $0\leq k \leq n-r-1$, \qquad  $n\geq r+1$.
\item{(b1)} There exists $\;\;n_0\geq r\;\;$ such that $\;\;\bigl<x^{n_0-r}Q_{n_0},
v\bigr>\not=\theta$. 

{\bf 2.-} We will say that $\bigl(Q_n\bigr)_{n\geq 0}$ is a   sequence
of regularly left quasi-orthogonal  matrix polynomials of order
$r$ with respect to $v \in \P^{{(m)}'}\setminus\{0\}$ if 
\item{(a)} $\bigl<x^k Q_n, v\bigr> = \theta, $\qquad $0\leq k \leq n-r-1$, \qquad  $n\geq r+1$.
\item{(b2)} There exists $\;\;n_0\geq r\;\;$ such that $\;\;\bigl<x^{n_0-r}Q_{n_0},
v\bigr>$ is non-singular.} 

  Notice that if $r=0$ the  previous definition is the classical orthogonality.  Furthermore, for
$m=0,$ both definitions yield the scalar quasi-orthogonality.

   Remember that  given   $ w\in \P^{{(m)}'}$ we will say that  $\;w=0\;$ if $\;\left<P,
w\right> = \theta$, $\;\;\forall P\in  \P^{(m)}$.

  As  we will see in the next proposition, when a sequence of quasi-orthogonal matrix polynomials
is already orthogonal with respect to another functional, the conditions (b1), (b2) become stronger.

\vskip 0.25cm

{\bf 3.2  Proposition.} {\sl Let $ u\in \P^{{(m)}'}$ be  quasi-definite and
let $\bigl(P_n\bigr)_{n\geq 0}$ be the co\-rres\-pon\-ding   sequence of left orthogonal matrix 
polynomials. Given $ v\in
\P^{{(m)}'}$, the sequence 
$\bigl(P_n\bigr)_{n\geq 0}$    is quasi-orthogonal of order $r$ with respect to $v$ if and only
if

\item {(a)}  $\bigl< x^k P_n, v\bigr> = \theta, $\qquad $0\leq k \leq n-r-1$, \qquad  $n\geq r+1$.

\item{(b'1)}   $\bigl<x^{n-r} P_n  , v\bigr> \not= \theta, $\qquad $\forall n\geq r$.}

  {\sl The sequence  $\bigl(P_n\bigr)_{n\geq 0}$ is regularly quasi-orthogonal of order $r$ with
respect to $v$ if and only if

 \item {(a)}  $\bigl< x^k P_n, v\bigr> = \theta, $\qquad $0\leq k \leq n-r-1$, \qquad  $n\geq r+1$.

\item{(b'2)}   $\bigl<P_n x^{n-r}, v\bigr>\;\; $ is non singular   $\forall n\geq r$.}

{ \bf{\sl  Proof.} } Let $n\geq r+1$ be verifying (b1) or (b2). Taking into account the recurrence
relation (1.1),  we obtain
$$
\gamma_{n+1}\left<x^{n-r} P_n, v\right> = - \beta_{n+1} \left<x^{n-r} P_{n+1}, v\right> -
\alpha_{n+1} \left<x^{n-r} P_{n+2},v \right> +$$
$$+ \left<x^{n-r+1} P_{n+1}, v\right> =  
\left<x^{n+1-r} P_{n+1}, v\right>
$$

{\noindent and} (b1) or (b2) holds for $n+1$. In a similar way,
$$
\gamma_n\left<x^{n-1-r} P_{n-1}, v\right> =   \left<x^{n-r} P_n  , v\right>
$$

{\noindent and} (b1) or (b2) is satisfied for $n-1$. $\diamond\diamond$

\vskip 0.25cm
  In what follows, the following result will be useful.

\vskip 0.25cm

{\bf 3.3  Lemma.} {\sl Let  $ u\in \P^{{(m)}'}$ be   quasi-definite   and  let 
$\left(P_n\right)_{n\geq 0}$ be  the co\-rres\-pon\-ding   sequence of  left  orthogonal matrix
polynomials.  Then, given 
$ w\in
\P^{{(m)}'},\;$     $w=0\;$ if and only if there exists $n_0$ such that 
$$
\left <x^k P_n, w\right> = \theta, \qquad 0\leq k \leq n, \qquad n\geq n_0.
$$}
{ \bf{\sl  Proof.} }  As a consequence of the recurrence
relation  (1.1),   
$$
<x^k P_{{n_0}-1}, w> =   
$$
$$
 = \gamma_{n_0}^{-1}[<x^{k+1} P_{n_0}, w> - \alpha_{n_0}<x^k  P_{{n_0}+1}, w> - \beta_{n_0}
<x^k  P_{n_0}, w>] = \theta
$$

\noindent{for} $0\leq k \leq n_0-1.$  So, the hypothesis is true for $n\geq  n_0-1,$ $0\leq k \leq
n$ too, and  as a consequence, for  all $n\geq 0,$  $k\geq 0.$ That means $<x^k, w>  = \theta$ for   all
$k\geq 0$ and then, $w=0.$ $\diamond\diamond$

\vskip 0.25cm
 Given a sequence of orthogonal matrix polynomials with respect to a functional $u,$ its
quasi-orthogonality  with respect to another functional $v$ can be cha\-rac\-te\-ri\-zed by a simple
relation beetwen $u$ and $v.$

\vskip 0.25cm

{\bf 3.4 Theorem.} {\sl Let $ u, v\in \P^{{(m)}'}$   such that $u$  is   quasi-definite
  and let 
$\left(P_n\right)_{n\geq 0}$ be  the corresponding     sequence of left orthogonal
matrix polynomials. Then,

 $\;\;$(i)  $\left(P_n\right)_{n\geq 0}$ is    quasi-orthogonal of order $r$ with respect to $v$ if
and only if there exists  
$A\in
\P^{(m)}$, with dg$A=r$,     and such that $v=uA$.

 $\;\;$(ii)  $\left(P_n\right)_{n\geq 0}$ is regularly   quasi-orthogonal of order $r$ with
respect to
$v$ if and only if there exists   
$A\in
\P^{(m)}$, with dg$A=r$,  and non-singular leading coefficient  such that $v=uA$.

  Moreover, in any case, the matrix polynomial $A$ is unique.}

{ \bf{\sl  Proof.} }

  We will do the proof just for (i) because for (ii) the arguments are   analogous.

$\underline{\Leftarrow}/$ We consider $\displaystyle A(x) = \sum_{j=0}^r a_j x^j$, with
$a_r\not=\theta$. Then,  for
$v=uA$
$$
\left<x^k P_n, v\right> = \sum_{j=0}^r\left< x^{k+j}P_n, u\right> a_j.
$$
 For $0\leq k\leq n-r-1$,  we have
$$
\left<x^{k+j} P_n, u\right> =  \theta,  \qquad \left(0\leq j\leq r\right)
$$
and, for  $k=n-r,$
$$
\left<x^k P_n, v\right> = \left<x^n P_n, u\right> a_r \not=\theta.  
$$

   So, $\;\;\left(P_n\right)_{n\geq 0}$ is regularly quasi-orthogonal of
order
$r$ with respect to $v$.

$\underline{\Rightarrow}/$  We  will prove that $\;\displaystyle 0=v-\sum_{j=0}^r
x^j u a_j\;\;$ has an unique solution in the unknowns $(a_j)_{j=0}^r$.

  For all $n\geq r$ fixed, the  system of $r+1$ equations
$$
\left\{\eqalign{&\left<x^{n-r}P_n, v\right> = \left<x^n P_n, u\right> a_r^{(n)},  \cr
&\cdots\cdots \cdots \cr
&\left<x^n P_n, v\right> = \left<x^n P_n, u\right> a_0^{(n)}+\dots+\left<x^{n+r} P_n, u\right>
a_r^{(n)},\cr}\right. \eqno (2.1)
$$

{\noindent has} a unique solution in the unknows  $(a_j^{(n)})_{j=0}^r$. Notice that the first
equation leads to $a_r^{(n)} \not=\theta.$   

  We  will prove that   
$(a_j^{(n)})_{j=0}^r$ are independent of $n$. Keeping in mind the recurrence relation (1.1), 
we have that
the  relation
$$
\gamma_{n+1} \left<x^kP_n, \hat v\right> + \beta_{n+1} \left<x^kP_{n+1}, \hat v\right>  +
\alpha_{n+1} \left<x^k P_{n+2}, \hat v\right> = \left<x^{k+1}P_{n+1}, \hat v\right> \eqno (2.2)
$$
 is  verified  $\forall k\geq 0$ and $\forall\hat v \in \P^{{(m)}'}$. The first equation of (2.1)
leads to
$$
\left\{\eqalign{&\left<x^{n-r}P_n, v\right> = \left<x^n P_n, u\right> a_r^{(n)}  \cr
                &\left<x^{n+1-r}P_{n+1}, v\right> = \left<x^{n+1} P_{n+1}, u\right> a_r^{(n+1)}, 
\cr}\right.\eqno (2.3)
$$
and,   from (2.2)
$$
\left\{\eqalign{&\gamma_{n+1}\left<x^{n-r}P_n, v\right> - \left<x^{n+1-r} P_{n+1},
v\right> = \theta, \cr
                &\gamma_{n+1}\left<x^n P_n, u\right> -  \left<x^{n+1} P_{n+1}, u\right>
 = \theta. \cr}\right.\eqno (2.4)
$$

   Equations (2.3) and (2.4)  imply   $ a_r^{(n)} = a_r^{(n+1)} $. 

  We supose $ a_l^{(n)} = a_l^{(n+1)} $
$(l=m+1,\dots,r, \quad m<r,\quad \forall n\geq r),$  and we will prove that $ a_m^{(n)} = a_m^{(n+1)}
$. Taking into account   (2.1) we have
$$
\left\{\eqalign{&\left<x^{n-m}P_n, v\right> = \sum_{j=m}^r\left<x^{n-m+j} P_n, u\right> a_j^{(n)} 
\cr
                &\left<x^{n-m+1}P_{n+1}, v\right> = \sum_{j=m}^r\left<x^{n+1-m+j} P_{n+1}, u\right>
a_j^{(n+1)} 
\cr}\right.\eqno (2.5)
$$
and from (2.2)
$$ 
\left\{\eqalign{\gamma_{n+1} \left<x^{n-m}P_n,v\right> + \beta_{n+1} &\left<x^{n-m}P_{n+1},v\right>
 +\alpha_{n+1} \left<x^{n-m} P_{n+2}, v\right> = \cr
 = &\left<x^{n-m+1}P_{n+1}, v\right>,\cr}\right.\eqno(2.6)
$$

$$ 
 \left\{\eqalign{\gamma_{n+1}\left<x^{n-m+j}P_n,u\right> + \beta_{n+1}
&\left<x^{n-m+j}P_{n+1},u\right> +  \alpha_{n+1} \left<x^{n-m+j} P_{n+2}, u\right> =\cr  
= &\left<x^{n+1-m+j}P_{n+1}, u\right>.\cr}\right.\eqno (2.7)
$$

 By substitution of (2.5) into (2.6) we obtain
$$
\gamma_{n+1}\sum_{j=m}^r \left<x^{n-m+j}P_n, u \right> a_j^{(n)} + \beta_{n+1}
\sum_{j=m+1}^r\left<x^{n-m+j}P_{n+1},u \right> a_j^{(n+1)}+$$
$$+ \alpha_{n+1}\sum_{j=m}^r\left<x^{n-m+j}P_{n+2}, u\right> a_j^{(n+1)} = \sum_{j=m+2}^r
\left<x^{n+1-m+j}P_{n+1}, u\right> a_j^{(n+2)}.
$$

 Now, keeping in mind the relation (2.7) and by application of the induction hypothesis, we have
$$
\gamma_{n+1}  \left<x^n P_n, u \right> a_m^{(n)} - 
   \left<x^{n+1} P_{n+1}, u \right> a_m^{(n+1)} = \theta.
$$

  Thus, by (2.4),  we conclude $ a_m^{(n)} = a_m^{(n+1)} $.

  Let us denote by $(a_j)_{j=0}^r$ the solutions of (2.1) for $n\geq r,$  that we have proved 
are independent  $n$.
If  $\displaystyle A(x) = \sum_{j=0}^r a_j x^j$ and $w= v - u A,$   we
have
$$
 \left<x^k P_n, w \right> =  \left<x^k P_n, v \right> - \sum_{j=0}^r  \left<x^{k+j} P_n, u \right>
a_j, \quad  0\leq k \leq n,\quad n\geq r.
$$

  When  $\;\;0\leq k \leq n-r-1$,   taking into account the orthogonality with respect to the
functional $u$ as well as the quasi-ortogonality with respect to the functional $v,$ we get  
$\left<x^kP_n, w\right>=\theta.\;\;$ If $n-r\leq k \leq n,$ then $\;\;\left<x^k P_n, w \right> =
\theta,\;\;$
because  $(a_j)_{j=0}^r$ are solutions of the equations (2.1).  So,
$w=0,$ according to  Proposition 2.2.

\vskip 0.25cm
%\vfill\eject
%\null\vskip 0.5cm
 
 To prove the uniqueness of the matrix polynomial $A,$ let us suppose    
  $A,B\in \P^{(m)}$   such that $ v=uA = uB.$ Then, from (i), it must be 
 $r=\;$dg$A =\; $dg$B$ since $(P_n)_{n\geq 0}$ can not be quasi-orthogonal of different order
with respect  to $v.$  Hence,
 $\displaystyle A(x)
=\sum_{j=0}^r a_j x^j$,  $\quad\displaystyle B(x) =\sum_{j=0}^r b_j x^j$ with $a_r, b_r\not=\theta$. 
From  
$$
\left<x^k P_n, v\right> = \sum_{j=0}^p \left<x^{k+j} P_n, u\right>a_j = \sum_{j=0}^q \left<x^{k+j}
P_n, u\right> b_j,
$$
and taking
$k=n-r,\dots,n,$ we have
$$
\left\{\eqalign{&\left<x^n P_n, u\right> \left(a_r - b_r\right) = \theta,\cr
                &\cdots \cdots \cdots \cr
                &\left< x^{n+r} P_n, u \right>\left(a_r - b_r\right) + \cdots + 
                 \left< x^n P_n, u \right>\left(a_0 - b_0\right) = \theta.\cr}\right.
$$

  Therefore,  $a_j  = b_j$, $(j=0,\cdots, r)$. $\diamond\diamond$

\vskip 0.5cm
%\vfill\eject

{\bf \S $\;\;\;$ 4 -   Semi-classical functionals. }

\vskip 0.25cm

 We consider the derivative operator on the space $\P^{{(m)}'}$ as the linear operator
$D:\P^{{(m)}'} \to \P^{{(m)}'}$ such that $\left<P, Du\right> = - \left< P', u\right>.$
From this and the definition of the right-product it is  straightforward to prove that
$\; D(u A)  =   (D u) A + u A',$  for all $u\in\P^{(m)'}$ and $A\in\P^{(m)}.$
   
\vskip 0.25cm

{\bf 4.1 Definition.}  {\sl Let $u\in \P^{{(m)}'}$ be  quasi-definite. We will say that
$u$ is   semi-classical   if there exist  $A, B\in \P^{(m)},$ with  det $A\not=0$, such that
it is verified the distributional equation $D\left(uA\right) = uB$.  We will also say 
that the
 co\-rres\-pon\-ding sequence of left  orthogonal
matrix polynomials 
$\left(P_n\right)_{n\geq 0}$ is
  semi-classical.}

\vskip 0.25cm

{\bf Remark. } Given  $u\in \P^{{(m)}'},$ the linear functionals $u_{i,j}:\P^{(m)}\to \C,$
$\;(i,j=1,\cdots,m)$ defined by $u_{i,j}(P) = u(P)_{i,j},$ $\;\forall P\in\P^{(m)}$ are called 
the components of $u.$  Then, $(uA)_{i,j} = \displaystyle\sum_{k=1}^m u_{i,k} A_{k,j}$  with 
the obvious 
definition for the multiplication of $u_{i,j}$ by a scalar polynomial.
Therefore,  det$A\not=0$ is the minimal requirement to ensure that  
the previous definition  will involve  to all the
components of the functional $u.$
\vskip 0.25cm
{\bf 4.2  Lemma.}   {\sl Let $u\in \P^{{(m)}'}$  such that  $D\left(uA\right) = uB.$
Then, for every   $C\in\P^{(m)}$,
$$
D\left(u A C \right) = u \left(  AC' + BC\right)
$$
}

{ \bf{\sl  Proof.} } It  is just a consequence of the rule for the derivation of the right-product. $\diamond\diamond$ 
\vskip 0.25cm
  The previous result implies that, for every $u\in\P^{(m)'},$ the set
$$
 {\Cal  A}_u = \{ A\in \P^{(m)} / D(uA) = uB, B\in \P^{(m)}\}
$$
is a right-ideal of $\P^{(m)}.$

  When dealing with   scalar semi-classical functionals we arrive in this way to an ideal of $\P,$ which is 
necessarily generated by a unique (up to non-trivial factors) polynomial. This generator is used 
to classify the 
scalar semi-classical functionals, ([11, 12]).

  Since in  the matricial case, a right-ideal of $\P^{(m)}$ is not necessarily principal, 
 we cannot use ${\Cal A}_u$
  for the classification of semi-classical functionals. However, if $u\in\P^{(m)'}$ is
semi-classical  then there exists $\alpha\in\P\setminus\{0\}$ (where $\P \equiv \P^{(1)}$)  such
that
$\alpha I\in {\Cal A}_u.$ To see this just notice that  if $A\in {\Cal A}_u$  with det$A\not=0,$
then (det$A) I\in{\Cal A}_u$ since  (det$A)I  = A A^+$ with 
 $A^+ =$adj$A \in\P^{(m)}.$ Therefore, for every semi-classical matrix functional $u\in\P^{(m)'},$ the set
$$
{\vartheta}_u=\{\alpha\in\P/ D(u\alpha I) = uB, B\in\P^{(m)}\}
$$
is a non-trivial ideal of $\P.$ We can use the \lq\lq essentially"  unique generator of this ideal to clasify 
the semi-classical  matrix functionals similarly to the scalar case.
\vskip 0.25cm
{\bf 4.3 Definition.} {\sl Let $u\in\P^{(m)'}$ be semi-classical and let $\alpha\in \P\setminus\{0\}$ be a 
polynomial with
 smallest degree such that  $D(u\alpha I) = u B,$ $B\in\P^{(m)}.$   Then, we say that $u$ is  of class 
$s= $max$\{p-2, q-1\},$ where $p=\;$dg$\alpha$  and $q=\;$dg$B.$}

\vskip 0.25cm
{\bf Remarks.}

(i)  The preceding discussion shows that this definition is well done: it always exists such a polynomial 
$\alpha$ and the definition of class does not depend on the choice of $\alpha.$

(ii) When $u\in\P^{(m)'}$ is semi-classical there exists $A\in{\Cal A}_u$ with det$A\not=0,$ but it 
is possible for the leading coefficient of $A$ to be singular.
However, there  always exists   $\widetilde{A}\in{\Cal A}_u$ with non-singular leading
coefficient. To see this just take 
$\widetilde{A}=$(det$A) I.$
 
\vskip 0.25cm
 The  following theorem gives the first characterization of semi-classical matrix functionals.   
\vskip 0.25cm
{\bf 4.4 Theorem.} {\sl Let  $\;\;u\in \P^{{(m)}'}$ be  quasi-definite and 
$\left(P_n\right)_{n\geq 0}$be  the co\-rres\-pon\-ding    sequence of left orthogonal matrix
polynomials. Then, the following
 statements are equivalent:}

{\sl (i)   $u$ is   semi-classical.}

{\sl (ii) There exists   $v\in\P^{{(m)}'}\setminus\{0\}$  such that 
$\;\;\left(P'_{n+1}\right)_{n\geq 0}\;$ is quasi-orthogonal with respect to $v$, and  
$\;\;\left(P_n\right)_{n\geq 0}\;\;$ is regularly   quasi-orthogonal with respect to $v.$ }

{ \bf{\sl  Proof.} } 

$(i)\;\Rightarrow \;(ii)$   Notice that  if $u$ semi-classical then   $D\left(u
\alpha  I\right) = u B,$ $\alpha\in\P\setminus\{0\}$ according to   previous comments.  Therefore,
$$
\left< \left(x^k P_n\right)', u\alpha I\right>   + \left< x^k P_n, u B\right>=\theta.
$$
   
   Let $\displaystyle \alpha(x)  = \sum_{j=0}^p a_j x^j$,   
 $\displaystyle B(x) = \sum_{j=0}^q b_j x^j$   with $ a_p\not=0 $ and  $ b_q\not=\theta$ 
whenever 
$B\not=\theta$.  With this notation the  above equation becomes
$$
 \sum_{j=0}^p\left<k x^{k-1+j}P_n, u \right> a_j +    
 \sum_{j=0}^q\left< x^{k+j}P_n, u\right> b_j  = - \left< x^k P'_n, u\alpha I\right>.
\eqno(3.1)
$$

 The left hand side is equal to zero matrix if $\;\;0\leq k \leq n-p\;\;$ and  $\;\;0\leq k\leq
n-(q+1).\;\;$ Let
$r=\;$max$\{p-1, q\}.$ Then, $\left<x^k P'_n, u\alpha I\right> = 0 $ if $ 0\leq k \leq n-r-1$.
For
$k=n-r,$     (3.1)    becomes
$$
(n-r)\left<x^n P_n, u\right> a_{r+1} +   \left<x^n P_n, u\right> b_r = - 
\left<x^{n-r} P'_n, u \alpha I\right>,
$$
and thus
$$
\left<x^{n-r} P'_n, u\alpha I\right> = -  \left<x^n P_n,  u\right> \left[(n - r) a_{r+1} +
b_r\right],
$$

{\noindent which} is  equal to zero matrix for  at   most one $n$. So, the sequence
$\left(Q_n\right)_{n\geq 0}$  with
$Q_n = P'_{n+1}$ is quasi-orthogonal   with respect to $u\alpha I,$ of order $r-1$.

  Obviously, since $\alpha I$ has non-singular leading coefficient, $\left(P_n\right)_{n\geq 0}$ is
regularly quasi-orthogonal   with respect to
$u\alpha I,$  of order $p$.

\vskip 0.5cm

$(ii)\; \Rightarrow \;(i)$   If $\left(P_n\right)_{n\geq 0}$  
 is  regularly quasi-orthogonal of order $p$ with respect to $v$, there exists $A\in \P^{(m)}$ 
with
dg$A=p$ and non-singular leading coefficient of $A,$ such that $v=uA$. Notice that if the leading
coefficient of $A$ is non-singular, then det$A\not=0$.

 Let $w= D(uA).$ It is,  $\left<P, w\right> = - \left<P', uA\right>$  and for $P=x^k P_n$ we
have 
$$
\left<x^k P_n, w \right> = -k \left<x^{k-1}P_n, u A\right> - \left<x^k P_n', u A\right>. 
$$
 
 From the quasi-orthogonality, the right hand side  vanishes if $\;0\leq k \leq n-p\;$ and
$\;0\leq k
\leq n-s-2,\;$ where $s$ is the order of quasi-orthogonality of $\left(P'_{n+1}\right)_{n\geq 0}$.
So,
$\;\left(P_n\right)_{n\geq 0}$ is quasi-orthogonal with respect to the functional
$w$ with order at  most max$\{p-1, s+1\}$. Thus, there exists  $B\in \P^{(m)}$ such that
dg$B\leq\;$max$\;\{p-1, s+1\}$  and $w = u B$. $\diamond\diamond$.

\vskip 0.25cm

{\bf Remark.}  Notice that if $\left(P'_{n+1}\right)_{n\geq 0}$ is quasi-orthogonal of order $s$
with respect to $u\alpha I,$ then  $D(u\alpha I) = u B$ where dg$\alpha=p,$  dg$B=q$ and this
implies that $s=\;$max$\;\{p-2, q-1\}.$

\vskip 0.5cm
%\vfill\eject
 
{\bf \S $\;\;\;$ 5 -   Structure relation and differo-differential equation. }

\vskip 0.25cm

{\bf 5.1 Theorem.}  (Structure  relation)  

 {\sl   Let  $\;\;u\in \P^{{(m)}'}$ be  quasi-definite and let $\left(P_n\right)_{n\geq 0}$  be
the
 associated    sequence of left orthogonal matrix  polynomials. Then, the following statements are
equivalent:

(i)    $u$ is   semi-classical.

(ii)   There exist  a polynomial $\alpha \in \P\setminus\{0\}$  with  dg$\alpha = p,$
   $s\in\N\cup\{0\},$ $\;\;s\geq p-2$ and $ \Theta_j^{(n)}
\in \C^{(m,m)}$  $(n\geq 0, -s\leq j\leq p),$   such that
$$
 \alpha(x) P'_{n+1}(x)  = \sum_{j=-s}^p \Theta_j^{(n)} P_{n+j}(x),
$$
where    $\Theta_{-s}^{(n)}  \not=\theta$  for some $n \geq s$ (we use the convention $P_k=\theta$ for
$k<0)$.}

{ \bf{\sl  Proof.} } 

$(i)\Rightarrow (ii)\;$  If $u$ is semi-classical then there exist $\alpha\in\P\setminus\{0\}$
 and $B\in\P^{(m)}$ such that 
$ D(u \alpha I) =  u B.$  Hence,   $\left(P'_{n+1}\right)_{n\geq 0}$ is quasi-orthogonal with
respect to $u \alpha I$ of order
$s=$max$\{p-2, q-1\},$  where $\;p=$dg$\alpha\;$ and  $\;q=$dg$B.\;$ This implies  
$$
\;\;\;<x^k P_{n+1}', u \alpha I> = \theta, 
 \qquad\;0\leq k\leq n-s-1,
$$ 
$$
<x^{n-s} P_{n+1}', u \alpha I> \not= \theta,\quad  for \quad some \quad n\geq s.
$$  
 
  Since    
$\alpha P'_{n+1}\in \P^{(m)}_{n+p}$   there exist  $\Theta_j^{(n)}$ verifying  $(ii)$.

\vskip 0.5cm

$(ii)\Rightarrow (i)\;$   Let $v=u\alpha I.$ Then  
$$
\left< x^k P'_{n+1 }, v\right> = \sum_{j=-s}^p\Theta_j^{(n)}\left< x^k P_{n+j}, u\right> = \theta,
\quad  0\leq k \leq n-s-1, \quad n\geq s,
$$
and, since $< x^{n-s} P_{n-s}, u>$ is non-singular,
$$
\left< x^{n-s} P'_{n+1 }, v\right> =  \Theta_{-s}^{(n)}\left< x^{n-s} P_{n-s}, u\right> \not=
\theta, \quad for \quad some  \quad n\geq s.
$$

 Thus, $\left(P'_{n+1}\right)_{n\geq 0}$ is quasi-orthogonal     of
order  $s$  with respect to $v.$ Obviously   $\left(P_n\right)_{n\geq 0}$ is 
regularly quasi-orthogonal of order $p$ with respect to $v$ and, so, $u$ is semi-classical.
$\diamond\diamond$

\vskip 0.5cm

{\bf 5.2 Theorem. }  (Differo-differential  equation)

 {\sl   Let  $\;\;u\in \P^{{(m)}'}$ be quasi-definite and let $\left(P_n\right)_{n\geq 0}$ be the
 corresponding   sequence of  left orthogonal matrix polynomials. Then, the following statements
are equivalent:

(i) $u$ is semi-classical.

(ii)  There  exist two polynomials 
  $\alpha,$ $\beta\in\P$ with $p=\;$dg$\alpha \geq 0,$ $q=\;$dg$\beta$ and matrices
$\Lambda_k^{(n)} \in \C^{(m,m)}$  $\left(n\geq 0, \; -s\leq k\leq s\right),$ 
such that 
$$
\alpha(x) P''_n(x)   + \beta(x) P'_n(x)   = \sum_{k=-s}^s \Lambda_k^{(n)} P_{n+k}(x),
$$
where $s\geq $max$\{p-2, q-1\}$ (we use  the convention $P_k=\theta$ for $k<0$).}  
% \vfill\eject

{ \bf{\sl  Proof.} }  

$(i)\Rightarrow (ii)$ $\;\;\;$ Let $u$ be semiclassical. By Theorem 5.1, there
 exist  a polynomial $a\in\P$, with dg$ a=p_1\geq 0$, a non-negative integer $s_1\geq p_1-2,$  and
matrices $\Theta_j^{(n)} \in \C^{(m,m)}$  $\;\left(n\geq 0,\;\;-s_1\leq j\leq p_1\right)$, such
that
$$
a(x) P'_n(x) = \sum_{j=-s_1}^{p_1} \Theta_j^{(n-1)} P_{n-1+j}(x),  \eqno(5.1)
$$
 
 Taking derivatives in (5.1), we obtain that
$$
a^2  P''_n + a  a'  P'_n  = \sum_{j=-s_1}^{p_1} \Theta_j^{(n-1)} a  P'_{n-1+j},
$$
  
  If we use (5.1) in the right hand side, it follows that
$$
a^2  P''_n  + a  a'  P'_n  = \sum_{j=-s_1}^{p_1} \Theta_j^{(n-1)} \sum_{k=-s_1}^{p_1}
\Theta_k^{(n-2+j)}  P_{n-2+j+k},
$$

  Let us denote $\alpha(x) = a^2(x),$ $\;\beta(x) = a(x) a'(x),$ $\; \Lambda_{j+k}^{(n)} = 
\Theta_j^{(n-1)} \Theta_k^{(n-2+j)},$ and $p=2p_1=$dg$\alpha,$  $q=2p_1-1=$dg$\beta,$
$s=$max$\{2p_1-2, 2s_1+2\}.$ Then,
$$
\alpha  P''_n  + \beta  P'_n  = \sum_{j=-s}^s \Lambda_j^{(n)} P_{n+j}, 
$$
with $s\geq$max$\{p-2, q-1\}.$

$ (ii)\Rightarrow (i) $ $\;\;\;$  Let $\left(\pi_n\right)_{n\geq 0}$ be the dual basis
of  $\left(P_n\right)_{n\geq 0},$ that is, $\pi_n = u P_n^* E_n^{-1}$, where $E_n:=<P_n,P_n>_u$.
Let  $\left(Q_n\right)_{n\geq 0}$ be  the basis of $\P^{(m)}$ given by $Q_n(x) =\displaystyle
{1\over n+1} P'_{n+1}$, and let $(\rho_n)_{n\geq 0}$ be the corresponding dual basis. Thus,

$$
D \rho_n = - (n+1) \pi_{n+1}, \qquad n\geq 0. \eqno(5.2)
$$

  We consider for every  $n\geq 0$ the linear functional $-D(\pi_n\alpha) + \pi_n\beta.$  
There exists $\left(\lambda_k^{(n)}\right)_{k\geq 0} \subseteq \C^{(m,m)}$ such that
$$
-D(\pi_n\alpha) + \pi_n\beta = \sum_{k=0}^\infty \rho_k \lambda_k^{(n)},
$$
where  from
$$
<Q_j, -D(\pi_n\alpha) + \pi_n\beta > =   \sum_{k=0}^\infty <Q_j, \rho_k> \lambda_k^{(n)} =
 \lambda_j^{(n)}
$$
holds for $j\geq 0$. Hence, our hypothesis implies that
$$
\lambda_j^{(n)} = {1\over n+1} \sum_{k=-s}^s  \Lambda_k^{(j+1)} \delta_{j+1+k, n}. \eqno(5.3)
$$

  For $n=0,$ (5.3) leads to $\lambda_j^{(0)} = 0$ if $j\geq s,$ and 
$$
-D(\pi_0\alpha) + \pi_0\beta = \sum_{k=0}^{s-1} \rho_k \lambda_k^{(0)},
$$
or
$$
-D(u\alpha) P_0^* E_0^{-1} + u\beta P_0^* E_0^{-1} = \sum_{k=0}^{s-1} \rho_k \lambda_k^{(0)}.
$$

  If we denote $\widetilde{\lambda}_k^{(0)} = \lambda_k^{(0)} E_0 (P_0^*)^{-1},$ then
$$
-D(u\alpha) + u \beta =  \sum_{k=0}^{s-1} \rho_k \widetilde{\lambda}_k^{(0)}\eqno(5.4)
$$

  In a similar way, for $n=1$ (5.3) gives
$$
-D(u\alpha) P_1^* E_1^{-1} + u\beta P_1^* E_1^{-1} = \sum_{k=0}^s \rho_k \lambda_k^{(1)}
$$
because $\lambda_j^{(1)}=0$ for $j\geq s+1.$  Thus,
$$
-D(u\alpha P_1^*) + u\beta P_1^* = \sum_{k=0}^s \rho_k \widetilde{\lambda}_k^{(1)} \eqno(5.5)
$$
with  $\widetilde{\lambda}_k^{(1)} = \lambda_k^{(1)} E_1.$

  Keeping in mind that $P_1(x) = M_1 x + M_2,$ with $M_1$ non-singular, (5.5) remains
$$
-D(u\alpha) P_1^* -u\alpha M_1^* + u\beta P_1^* = \sum_{k=0}^s \rho_k \widetilde{\lambda}_k^{(1)}, 
$$
and,  using (5.4),
$$
- u \alpha M_1^* = \sum_{k=0}^s \rho_k \widetilde{\lambda}_k^{(1)} - \left(
\sum_{k=0}^{s-1} \rho_k \widetilde{\lambda}_k^{(0)}\right) P_1^*.
$$

  Taking derivatives and applying (5.2) we obtain
$$
- D(u\alpha) M_1^* = - \sum_{k=0}^s (k+1) \pi_{k+1} \widetilde{\lambda}_k^{(1)}  +
\left(\sum_{k=0}^{s-1} (k+1) \pi_{k+1}
  \widetilde{\lambda}_k^{(0)}\right) P_1^* -
\left(\sum_{k=0}^{s-1} \rho_k \widetilde{\lambda}_k^{(0)}\right) M_1^* =
$$
$$
= -\sum_{k=0}^s (k+1) u P_{k+1}^* E_{k+1}^{-1} \widetilde{\lambda}_k^{(1)}+
\left(\sum_{k=0}^{s-1} (k+1) u P_{k+1}^* E_{k+1}^{-1} \widetilde{\lambda}_k^{(0)}\right) P_1^* -
\left(\sum_{k=0}^{s-1} \rho_k \widetilde{\lambda}_k^{(0)}\right) M_1^*.
$$

   As a  consequence, the polynomial $\Psi\in \P^{(m)}$ given by
$$
\Psi(x):= \left[\sum_{k=0}^s (k+1)  P_{k+1}^*(x) E_{k+1}^{-1} \widetilde{\lambda}_k^{(1)} -
\left(\sum_{k=0}^{s-1} (k+1) P_{k+1}^*(x) E_{k+1}^{-1}  \widetilde{\lambda}_k^{(0)}\right)
P_1^*(x)\right]\left(M_1^*\right)^{-1},
$$
whose degree is  not greater  than $s+1,$ verifies that
$$
D(u\alpha) = u \Psi + \sum_{k=0}^{s-1} \rho_k \widetilde{\lambda}_k^{(0)}.
$$

  From this and  (5.4) we  get 
$$
D(u\alpha) = u B
$$
with $B=\displaystyle {\Psi + \beta I \over  2},$  that is, $u$ is semi-classical.$\diamond\diamond$

\vskip 0.5cm

{\bf    Remark 1. } The    distributional equation  $D(u\alpha) = u \beta $ for  a  scalar
func\-tio\-nals $u$  implies that

$$
<x^k(\alpha p_n''(x) + \beta p_n'(x)), u> =  
$$
$$
= <(x^k p_n'(x))', u \alpha> - k<x^{k-1} p_n'(x), u
\alpha> + < x^k p_n'(x), u\beta> =
$$
$$
= - k <x^{k-1} p_n'(x), u\alpha >,
$$
wich vanishes  for $k= 0,\dots, n-s-1.$ So,

$$
\alpha p_n''(x) + \beta p_n'(x) = \sum_{k=-s}^s \lambda_{nk} p_{n+k}(x)
$$

  In the matricial case, the equality    $D(u\alpha I ) = u B $ does not imply $\;\;   <x^k P_n'(x),$
$ \;\; u B(x)>
 = <x^k B(x) P_n'(x), u>$
 for the non-conmutativity, and more generaly, there not  exists a polynomial $B(x;n)$ such
that  $\;\;<x^kP_n'(x), u B(x)>\;\;  = \;\; $
$\;\; < x^k B(x;n) P_n'(x), u>,$ neither.

  In this situation it  is not   possible to obtain a  differo-differential equation with  the
polynomial $B$ explicitly. 

\vskip 0.5cm

{\bf 5.2  Remark 2. }  Obviuosly, the differo-differential equation   given by the Theorem  5.2
is  not  unique because we can modify it by adding any structure relation.

\vskip 0.5cm

{\bf \S $\;\;\;$ 6 -   Some examples. }

\vskip 0.25cm
   
  We will illustrate the preceding results with some examples. 
\vskip 0.25cm

  {\bf 1.} In the first one, we consider the scalar
Laguerre weight $w= x e^{-x}$ that verifies the  Pearson type equation
$$
( w x )' = w ( 2 - x).
$$

 Now, we define the matrix weight function
$$
u:=w \pmatrix  
1&1\\1&1+x^2
\endpmatrix
$$
that also verifies a Pearson type equation. In fact
$$
D(u x I) = u B(x)
$$
where $B(x)= \pmatrix 2-x&-2\\0&4-x \endpmatrix$.

 Moreover, there exists $T\in GL(\R^2)$ such  that
$$
\hat u = T u T^t = \pmatrix  x e^{-x}& 0\\ 0&x^3 e^{-x} \endpmatrix, \qquad  T =
 \pmatrix 1&0\\ -1&1\endpmatrix
$$
i.e., $u$ is congruent with a diagonal  weight  with classical scalar weights in
its diagonal and then, the functional $u$ is positive definite.

  Let $l_n^{(a)}$ be the monic Laguerre polynomials asociated to the weight  $x^a e^{-x},$ $
a > -1.$  Then,
$$
L_n(x) =  \pmatrix l_n^{(1)}(x)&0 \\ 0&l_n^{(3)}(x) \endpmatrix
$$
constitutes a sequence of matrix O P   related to the diagonal weight $\hat u.$ Obviously 
$\hat u$ satisfies  a Pearson  equation,
$$
D( \hat u x I ) = \hat u \hat B(x), \qquad  \hat B(x) = ({T^t})^{-1} B(x) T^t =  \pmatrix 2-x
&0\\ 0&4-x \endpmatrix,
$$
and a sequence of matrix O P associated to  $u$ is given by  
$$
P_n(x) = L_n(x) T, \qquad  n\geq 0
$$

  It is easy to obtain a estructure relation
$$
x P_{n+1}'(x) = (n+1) P_{n+1} + (n+1)  \pmatrix n+2&0\\ 0&n+4\endpmatrix P_n(x),
$$
as well as, for all $n\geq 2,$ a differo-differential equation,
$$
x^2 P_n''(x) + x P_n'(x) = n (n+1) P_n(x) +
$$
$$
+ n^2  \pmatrix n+1&0\\ 0&n+3 \endpmatrix P_{n-1}(x) +
n (n-1)   \pmatrix n(n+1)&0\\ 0&(n+2)(n+3) \endpmatrix P_{n-2}(x).
$$

  Moreover, the polynomials $l_n^{(a)}$  satisfy  the differential equation,
$$
x l_n^{(a)''}(x) + ( x - a - 1 )  l_n^{(a)'}(x) = n l_n^{(a)}(x),
$$
and so, we can obtain for the polynomials  $P_n$ the following differential equation:
$$
x P_n''(x) +  \pmatrix x-2&0\\ 0&x-4 \endpmatrix P_n'(x) = n P_n(x).
$$

 For the monic polynomial     $\widetilde{P}_n = T^{-1} P_n$ we have  
$$
x\widetilde{P}_n''(x) + T^{-1}  \pmatrix x-2&0\\ 0&x-4 \endpmatrix T \widetilde{P}_n'(x) = n
\widetilde{P}_n(x),
$$
and then
$$
x\widetilde{P}_n''(x)  + B^t(x) \widetilde{P}_n'(x) = n \widetilde{P}_n(x)
$$

  The polynomials $\widetilde{P}_n$ can be obtained by means of a Rodrigues type formula

$$
\widetilde{P}_n(x) = {1\over n!} \left(T^t T\right)^{-1} u^{-1} {d^n\left(x^n u\right)\over d x^n}
\left(T^t T\right).
$$

\vskip 0.5cm
   
 {\bf 2.}  Now,  we consider the 
      the matrix  weight
$
u:= e^{ - x^2} R(x)
$
where 
$
R(x) =  \pmatrix 0&1  \\  1& x \endpmatrix
$. So,  $R(x)$  is not   positive definite.

 Then, $u$ verifies the distributional equation   
$$
u' = e^{ - x^2} R'(x) - 2 x e^{ - x^2} R(x) = u( -2x I + R^{-1}(x)R'(x)) = u B(x), \eqno(6.1)
$$
where 
$
B(x) =  \pmatrix -2x&1  \\  0&-2 x \endpmatrix
$. Moreover, $u$  has not a congruent diagonal.
  
 Now, let
$$
P_n(x) = (-2)^{(-n)} e^{x^2} S^{-1}(x){d^n(e^{-x^2} S(x)) \over dx^n}, \quad n\geq 0,\eqno(6.2)
$$
where 
$
S(x) = 
 \pmatrix
1+x&1  \\  x&1  
\endpmatrix 
$.

  By   derivation  in (6.2)
$$
P_{n+1}(x) = -{1\over 2}[-2xI + S^{-1}(x) S'(x)] P_n(x) - {1\over 2} P_n'(x).
$$

  Then,
$$
P_{n+1}(x) =  -{1\over 2} [B^t(x) P_n(x) + P_n'(x)], \eqno(6.3)
$$
and the leading coefficient of $P_n$ is $I$ because $P_0(x) = I.$

 Applying  the Leibnitz rule for the derivatives in (6.2), we have
$$
P_{n+2}(x) = (-2)^{-n-2} e^{x^2} S^{-1}(x) {d^{n+1}(-2x e^{-x^2} S(x) + e^{-x^2} S'(x))\over
dx^{n+1}} =
$$
$$
=  (-2)^{-n-1} x e^{x^2} S^{-1}(x) {d^{n+1} (e^{-x^2} S(x)) \over dx^{n+1}} +
(n+1) (-2)^{-n-1} e^{x^2} S^{-1}(x) {d^n (e^{-x^2} S(x)) \over dx^n} +
$$
$$
+  (-2)^{-n-2} e^{x^2} S^{-1}(x) S'(x) S^{-1}(x) {d^{n+1} (e^{-x^2} S(x)) \over dx^{n+1}}.
$$

 So, we obtain the recurrence relation,
$$
P_{n+2}(x) =  -{1\over 2} [B^t(x) P_{n+1}(x) + (n+1) P_n(x)]. \eqno(6.4)
$$

 As a consequence, $(P_n)$ is a  sequence of matrix orthogonal polynomials with respect to  
certain  non-singular functional  (non-positive definite because $-{1\over 2}(n+1) I$ is a non-
positive definite matrix).

  Moreover,  from (6.3) and (6.4) we can obtain a structure relation
$$
P_n(x) = n P_{n-1}(x), \eqno(6.5)
$$
and by  derivation, the differo-differential equation
$$
P_n''(x) = n(n-1) P_{n-2}(x). \eqno(6.6)
$$

  So, $(P_n)$ is a semiclassical matrix orthogonal polynomial sequence.

  Now, from (6.3) and (6.5) it follows the  differential equation,
$$
P_n''(x) + B^t(x) P_n'(x) = - 2n P_n(x). \eqno(6.7)
$$

  Indeed, $(P_n)$ is an  sequence of matrix orthogonal polynomials with respect to the   matricial
 weight $u.$  In fact, 
$$
\int_{-\infty}^{\infty} P_n'(x)   e^{-x^2} R(x) P_m^{\prime t }(x) dx =
$$
$$
= \int_{-\infty}^{\infty} [ P_n''(x) + (-2x I + R'(x) R^{-1}(x)) P_n'(x)]  e^{-x^2} R(x) P_m^t(x)
dx =
$$
$$
= - 2n \int_{-\infty}^{\infty} P_n(x)   e^{-x^2} R(x) P_m^t(x) dx, 
$$
the last equality from equation (6.7).  In the same way, we obtain for the above integral that
$$
\int_{-\infty}^{\infty} P_n'(x)   e^{-x^2} R(x) P_m^{\prime t }(x) dx =
- 2m \int_{-\infty}^{\infty} P_n(x)  e^{-x^2} R(x) P_m^t(x) dx. 
$$

  So, when $n\not=m$ we have the orthogonality for $(P_n)$ with respect to $u$. 

  As a consequence we
have a quasi definite functional such that:

(i) It satisfies the distributional equation $D(u\alpha I) = u B.$

(ii) The corresponding sequence of matrix orthogonal polynomials satisfies a Rodrigues type
formula.

(iii) It verifies  the structure relation (6.5), the differo-differential equation (6.6) and a
differential equation (6.7). Keeping in mind  Remark 1 to Theorem 5.2, notice that we have (6.7)
because   $u B = (u B)^t = B^t u  $   and   $ B^t P_n = P_n B^t.  $

\vskip 0.75cm
   
 {\bf 3.}  In the last   example, we consider the Jacobi scalar weight $w=1-x^2$ that satisfies
the distributional equation
$$
\left(w(1-x^2)\right)' = - w \cdot 4x, \qquad x\in[-1,1].
$$

  We define the matricial weight function,
$$
u:= w \pmatrix 1&x\\ x&2-x^2\endpmatrix, \qquad x\in[-1,1],
$$
that verifies the Pearson type equation

$$
D\left(u(1-x^2) I\right) = u B(x),
$$
where

$$
B(x) =  \pmatrix -9x& 2+x^2\\ 1&-11x \endpmatrix.
$$

  Notice that $u$ is positive definite and  then,  there exists a sequence of matrix orthogonal
polynomials $(P_n)$ related to $u.$ Moreover, $u$ not is congruent with a dia\-go\-nal weight. So,
the sequence $(P_n)$ satisfies a  structure relation with $s=1,$ $p=2$ and a  differo-differential
equation of  second order in the following way:
$$
(1-x^2)^2 P_n''(x) - x(1-2x) P_n'(x) = \sum_{k=-4}^2 \Lambda_{nk} P_{n+k}.
$$

%\vskip 0.75cm  

\vfill\eject 
\null\vskip 0.5cm
{\bf Acknowledgements}.- This research was   supported by 
 Direcci\'on General de Ense\~nanza Superior of Spain (DGES). Project PB98-1615.

   The authors are very grateful to Professor {\it  Francisco Marcell\'an} for his  remarks and
useful suggestions.

\vskip 0.75cm
  
\overfullrule = 0pt

{\bf References}

\vskip 0.5cm

 [1] M.  Alfaro,  A.  Branquinho,   F.  Marcell\'an, and  J. Petronilho.   
  {\rm A generalization of a theorem of S. Bochner.} Publicaciones del Seminario Matem\'atico
Garc\'{\i}a de Galdeano. Serie II. Secci\'on 1, n\'umero 11. (1992). 
\vskip 0.25cm

 [2] S. Bonan, D.S. Lubinski, and P. Nevai.  {\rm Orthogonal polynomials and their derivatives.}  
SIAM {\it  J. Math. Anal. }{\bf 18}  (1987) 1163-1175.
\vskip 0.25cm

 [3] T.S. Chihara.    
  {\it \lq\lq An introduction to orthogonal polynomials"}
  {\rm Gordon and Breach, New York, 1978.}  
\vskip 0.25cm

 [4] A.J.  Dur\'an.  
{\rm On orthogonal  polynomials  with respect to a positive definite matrix of measures.}
{\it Can. J. Math.} {\bf 47} (1995)  88-112.
\vskip 0.25cm

  [5]   A.J.  Dur\'an.  {\rm    Matrix inner product having a matrix symmetric second order 
differential  equation.} {\it Rocky Mountain Journal of Mathematics} {\bf 27} (1997) 585-600.
\vskip 0.25cm

  [6]   A.J. Dur\'an, W.  Van Assche. 
{\rm Orthogonal matrix polynomials and higher-order recurrence relations.}
{\it Linear Alg. Appl.} {\bf 219} (1995) 261-280.
\vskip 0.25cm

  [7]   G. Freud.  {\it \lq\lq  Orthogonal polynomials"}
  {\rm Pergamon Press, Oxford, 1971.} 
 \vskip 0.25cm 
  
 [8]   Y.L. Geronimus.   {\it  \lq\lq  Orthogonal polynomials"}
  {\rm Consultants Bureau, New York, 1961.}    
 \vskip 0.25cm 

 [9]  E. Hendriksen, H. van Rossum.  {\rm  Semi-classical orthogonal polynomials. }  
C.Brezinski et al. Eds., 
{\it Lecture Notes in Math.}{\bf 1171} (Springer, Berlin, 1985)   354-361.
\vskip 0.25cm 

 [10]  F. Marcell\'an, H.O.  Yakhlef.  {\rm  Recent trends on analytic properties of matrix
orthonormal polynomials. }  {\it  Electronic Transactions in Numerical Analysis.} To appear.
 
\vskip 0.25cm

 [11]   P. Maroni.  {\rm Variations around classical orthogonal polynomials. Connected problems. } 
{\it J. Comput. Appl. Math. }{\bf 48} (1-2) (1993) 133-155.
\vskip 0.25cm

 [12]  P. Maroni.  {\rm Une th\'eorie alg\'ebrique des polyn$\hat o$mes orthogonaux. Application
aux polynomes orthogonaux semiclassiques. C. Brezinski et al. Eds. }  {\it Orthogonal Polynomials
and Their Applications, IMACS Ann.Comput.Appl. Math. }{\bf 9} (1991) 95-130.
\vskip 0.25cm

 [13]  G. Sansigre. 
 {\it  Polinomios ortogonales matriciales y matrices bloques. }
 	{\rm Doctoral Dissertation
   Universidad de Zaragoza, (1992). In Spanish.}
\vskip 0.25cm  
 [14]  J.A. Shohat.   {\rm A differential equation for orthogonal polynomials. } 
{\it Duke. Math. Journal. }{\bf 5}  (1939) 401-407.

\bye
\end